# Depth of planning the state of a dynamic discrete system by autocorrelation function

Sergei Masaev
Institute of Oil and Gas
Siberian Federal University
Krasnoyarsk, Russia
Smasaev@sfu-kras.ru

*Abstract*— **The production system (multidimensional object) is considered as a dynamic system with discrete time. Formalized: space (state of the object, control actions, goals, observed values, analytical estimates). Analytical estimates of the state of a dynamic system are formed through the autocorrelation function. The autocorrelation function is calculated with the regulator setting the length of the analyzed time series (analysis depth). A digital copy of the production system is created, characterized by 1.2 million parameters. Modeling the activities of the production system is performed in the author's complex of programs. In total, twenty-eight controller states are calculated to analyze the effect of repeating parameters affecting the activity of the production system. The simulation shows the cyclical dynamics of changes in the autocorrelation function. Formalization of the production system is carried out, which allows you to move on to other methods of analysis of the production system: Kalman filter, neural network forecast, recurrence equation, balances.**

## I. INTRODUCTION

There are complex multidimensional objects, including the regulation of human activities, for example, production systems. For such objects in control theory [1-8], the following were considered: a mathematical model of optimal controlled processes [1, 2], sufficient optimality conditions, numerical optimization methods for multi-step processes with discrete control, solutions of some particular classes of discrete programming problems [9, 10].

Actual issues for consideration are: the control loop of multidimensional production systems, decentralized control of multidimensional production systems.

## II. FORMULATION OF THE PROBLEM

Without a feedback control loop, it is impossible to imagine the successful development of a production system in a dynamically changing environment. Traditionally, when managing a production system, a person strives to develop a system of measures for a different plan horizon.

The aim of this work is to supplement the method of integral indicators with values of the autocorrelation function to determine the planning horizon of the state of the production system under the influence of the external environment.

Achieving the goal depends on the tasks:

- formalize the production system as a dynamic system with discrete time;
- develop a digital copy of a dynamic system (production system);
- calculate the autocorrelation function of a set that displays the state of a digital dynamic system (production system).

Along with the classical methods of control theory, the original method for analyzing the activity of the production system turned out to be the method of integral indicators, first set out in 2009 to assess the state and managerial decisions in production systems.

## III. RESEARCH METHODOLOGY

We have a system $S$

$$S=\{T,X,U,Y,J,V,\phi,\psi,h,\xi\},$$

where

$T=\{t/t=0,1,2,...\}$ – discrete set of moments of time (moments of accounting system);

$X$ – system phase area, $x(t)=\left[x^1(t),x^2(t),...,x^n(t)\right]^T\in X$ – $n$ - vector of phase variables that determine the state of the system, and $x^*(t)=\left[x^{*1}(t),x^{*2}(t),...,x^{*n}(t)\right]^T\in X$ – $n$ - vector of phase variables that determine the requirements for the state of the system when allocating resources by experts or by the Bellman method. The phase variables of the system $x^n(t)$ are the values of the input or expenditure of resources.

$U$ – system control area, $u(t)=\left[u^1(t),u^2(t),...,u^m(t)\right]^T\in U$ - $m$ - control vector;

$Y$ – space of observed values with operational accounting of the system, $y(t)=\left[y^1(t),y^2(t),...,y^k(t)\right]^T\in Y$ - $k$ - vector of operational accounting system parameters;

$J$ – space control purposes, $J(t)=\left[j^1(t), j^2(t),..., j^c(t)\right]^T \in J$ - $c$ - system operation target vector during system operational control;

$V=V^+ \cup V^-$ – analytical space, $V(t)=\left[v^1(t), v^2(t),..., v^s(t)\right]^T \in V$ - $s$ - vector of analytical estimates;

$\phi: T\times X\times U \to X$ - transition function of the system can be represented as:

$$x(t+1)=\phi(x_0(t), x(t), u(t)),$$

where $x_0=x(0)$.

$\psi: T\times X \to Y$ - observation function defining parameters available for observation:

$$y(t)=\psi(x(t)),$$

$v^-: T\times X \to V$ - function of the analysis of the system in past periods of time:

$$v^-(t)=\xi^-(x(t-1), x(t-2),..., x(t-k)),$$

$v^+: T\times X \to V$ - function for analyzing system activities at future times:

$$v^+(t)=\xi^+(x_0(t+1), x_0(t+2),..., x_0(t+k))$$

Control $\phi(x_0, x(t), u(t))$ can be specified by a vector equation

$$x(t+1)=A(t)x(t)+B(t)u(t)+\xi(t), \qquad (1)$$

where

$A(t)$ - $n\times n$ is the matrix that determines the internal structure of the system;

$B(t)$ - $n\times m$ matrix defining the control structure;

$x(0)=x_0$, $\xi(t)$ - noise, random interference, environmental influences.

The method of integral indicators for analytical evaluation uses the parameter values $x(t)$ from $k$ previous measures to calculate the cross-correlation of parameters. The parameter $k$ is the depth of analysis, then we can form a matrix $X_k(t)$. In our case, instead of integral indicators, we will use the autocorrelation function of the parameters $x(t)$. Next, we center and normalize its elements $\overset{o}{X}_k(t)$ and we can calculate the autocorrelation function

$$W_K(t)=\frac{1}{k-1}\overset{o}{X}{}_k^T(t)\overset{o}{X}_k(t)=\left\|w_{i,i}^k(t)\right\|, \qquad (2)$$

$$w_{i,i}^k(t)=\frac{1}{k-1}\sum_{l=1}^{k}\overset{\circ}{x}{}^i(t-l)\,\overset{\circ}{x}{}^i((t+c)-l), \qquad (1)$$

where

$k$ - depth of analysis (planning horizon),

$w_{i,i}^k$ - function autocorrelation, $i=1,2,3,...,n$.

$c$ - selectable billing period number, $c=1,2,3,...,t$.

The first period of the matrix $W_K(t)$ is taken as the base of relative values $w_{i,i}^k$ in subsequent periods $t$. There is an opportunity to track the dynamics of changes in the autocorrelation function under the influence of repeating factors in the function $i$. For $2\le k\le K$ we get $K-1$ matrices $W_K(t)$ for the considered period of time $T$:

$$W_K(T)=\begin{bmatrix} w_{1,1}^2(1) & w_{1,1}^2(2) & \cdots & w_{1,1}^2(t) \\ w_{2,2}^3(1) & w_{2,2}^3(2) & \cdots & w_{2,2}^3(t) \\ \vdots & \vdots & \ddots & \vdots \\ w_{i,i}^k(1) & w_{i,i}^k(2) & \cdots & w_{i,i}^k(t) \end{bmatrix}, \qquad (2)$$

Where $K=T/2$.

We form a vector for each matrix from the sum $w_{i,i}^k$ over $t$ periods:

$$W_k(t)=\sum_{i=1}^{n} w_{i,i}^k(t) \qquad (3)$$

Following the method of integral indicators, the value $W_k(t)$ is decomposed (8) into the position value $W_k^+(t)$ and negative value $W_k^-(t)$, which form the integral indicators of the autocorrelation of the function over the base period:

$$W_k^+(t)=\sum_{i=1}^{n} w_{i,i}^k(t):(w_{i,i}^k(t)>0) \qquad (4)$$

$$W_k^-(t)=\sum_{i=1}^{n} w_{i,i}^k(t):(w_{i,i}^k(t)<0) \qquad (5)$$

$$W_k(t)=W_k^+(t)+W_k^-(t) \qquad (6)$$

From the obtained vectors we form a matrix:

$$W_K(T)=\begin{bmatrix} W_2(1) & W_2(2) & \cdots & W_2(t) \\ W_3(1) & W_3(2) & \cdots & W_3(t) \\ \vdots & \vdots & \ddots & \vdots \\ W_k(1) & W_k(2) & \cdots & W_k(t) \end{bmatrix} \qquad (7)$$

Summarize the elements $W_k$ of the column $t$:

$$W_k = \sum_{k=1}^{n} W_k(t) \tag{8}$$

Accordingly, the vector $W$ is transposed and equal to:

$$W^T = [W_1, W_2, ..., W_k] \tag{9}$$

## IV. Description of the production system

The production system of deep wood processing in the Krasnoyarsk Territory is considered [11]. The volume of processed raw materials is 1,000 thousand cubic meters are harvested of round wood in the North Yenisei region. Harvested raw materials are delivered on barges along the Yenisei River during the shipping period from June to September. From the harvested round wood, the following products are produced: floorboards, glued beams, euro-lining and etc. Matrix *A* is determined by the interconnection of parameters: harvesting complexes, timber trucks, wood processing line, pellet equipment, boilers and maintenance equipment for the main production. The production process is carried out with the mandatory completion of the procurement phase of the formation of a sufficient amount of raw materials and semi-finished products in warehouses. By the 25th period, equipment is purchased to increase the procurement, processing of raw materials by 2 times and maintain production at the required technical level. In the second year of work, there is a 2-fold increase in production. Before increasing production, measures are being taken to increase the efficiency of production management processes. One of the important properties of the production system is the decentralized development of plans and the use by decision centers of planning methods that are not related to control theory. The control matrix *B* is defined by the convenient syntax of decision centers for control in the production system [12-17].

## V. Research results

Modeling of the production system is performed in the author's software package.

A digital copy of the object has the following characteristics: *X*=5,675,805, *k*=adjustable, *N*=1.2 million, *T*= 60, *t*=1, $\xi$ - the environmental impact in the form of fines and unpredictable additional resource costs [18].

In Fig. 1 the autocorrelation of functions (5) was calculated according to which factors affecting the activity of the production system are seasonal in nature.

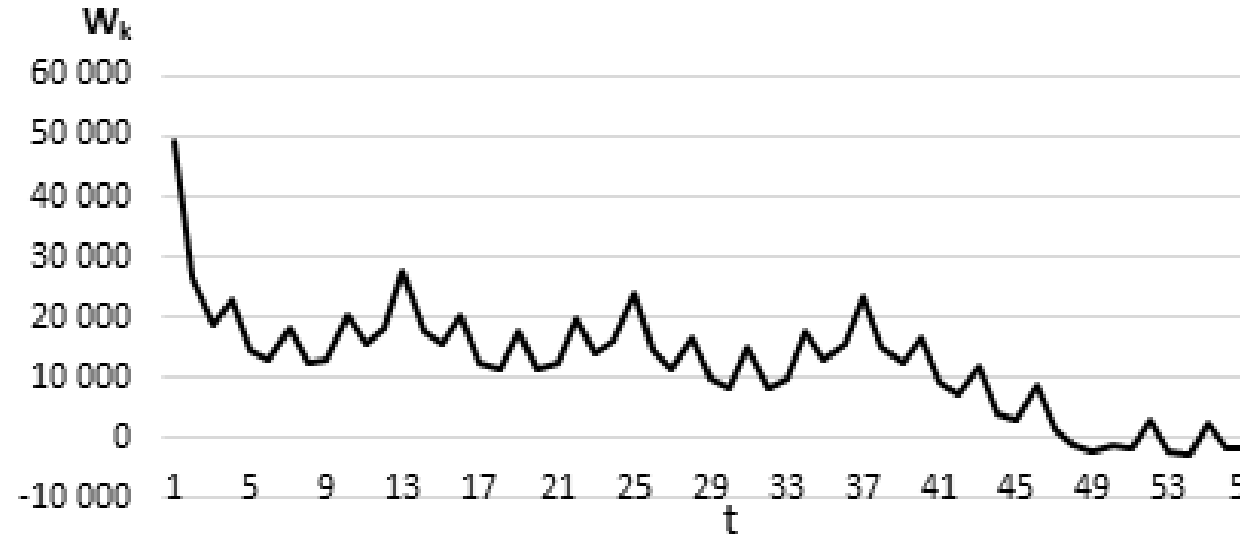

Fig. 1. The dynamics of the integral indicator $W_k$ for 5 years

By setting the magnitude of the controller $k$, an increase in the value of fluctuations of the functions $x^i$ is recorded, which makes it possible to determine the necessary depth (planning horizon) of the analysis of the studied system. This parameter $W_K$ is calculated (10). In Fig. 2 the calculation results are presented.

In Fig. 2 $W_K(T)$ formed from $W_k(t)$ for each value $k$ from 0 to 30. In Fig. 3 for clarity, we display $W_k(t)$ for each $k$. It can be seen from it that it is impossible to get by with a large horizon of event analysis, since there are fluctuations (factors) in autocorrelation that occur only once in two periods, which appear with the values of the controller $k=2$ or $k=6$.

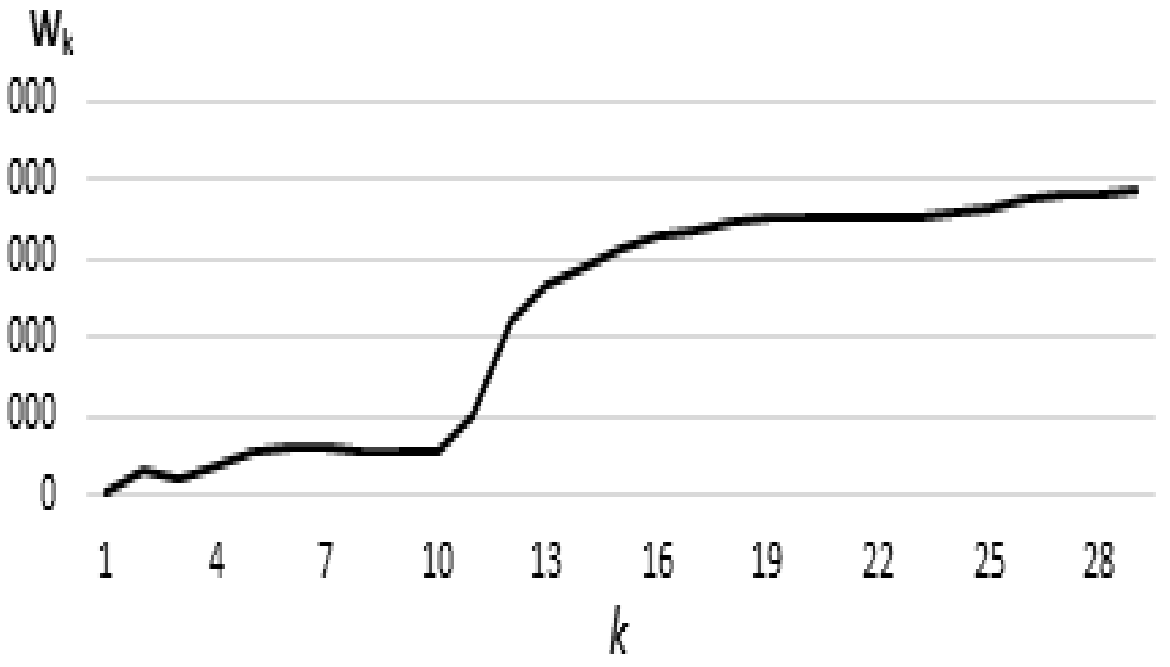

Fig. 2. The dynamics of the integral indicator $W_k$ depending on the values $k$

In Fig. 3, the maximum values $W_k(t)$ are visible in periods 13, 25 and 37. Through analysis of the functions $x_i$, it was found that the fluctuations are caused by an external factor of 6% on products and cost. In periods 3, 6, 9, 12, 15, 18, 21, 23, 26, 29, 32, 35, 38, 41, 43 and 46, self-oscillations of functions from fines operating from the external environment were recorded.

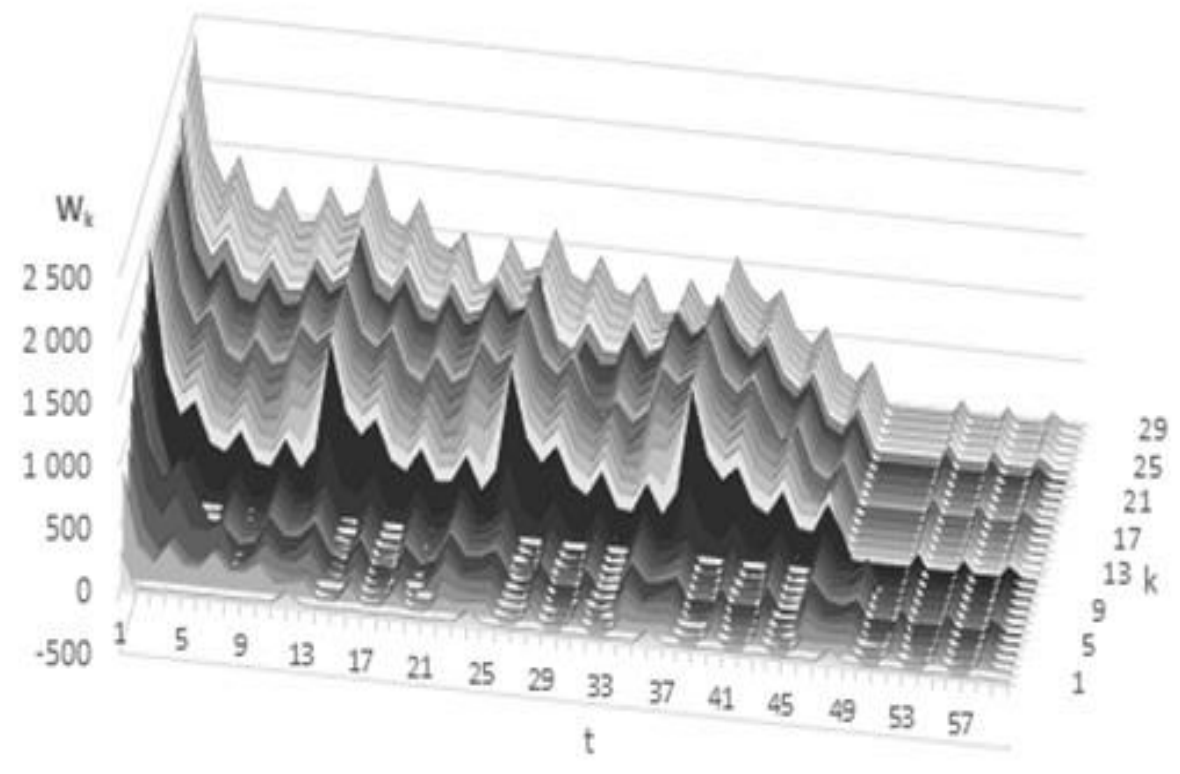

Fig. 3. Dynamics of the integral indicator $W_k$ for periods $t$ depending on the values $k$

Following the method of integral indicators $W_k(t)$ are calculated from positive $W_k^+(t)$ and negative $W_k^-(t)$ values, (dashed line) and (line) according to (6) and (7), respectively. In Fig. 4 the result is shown.

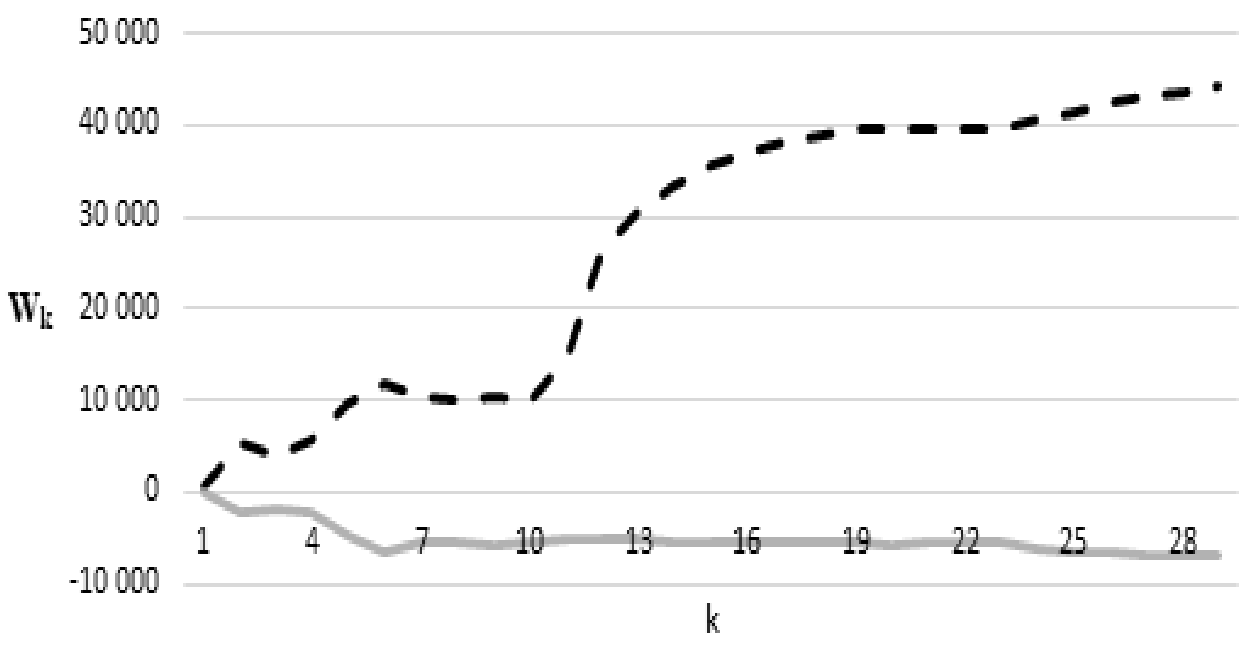

Fig. 4. Dynamics of positive and negative integral indicators $W_K$ depending on the values $k$

In Fig. 4 the dynamics $W_k^+(t)$ repeats the dynamics $W_k$, therefore, the choice of the regulator value $k$ remains the same with

- analysis of the short-term horizon (operational level of control) $k = 2$;
- analysis of the medium-term horizon (tactical level of control) $k = 6$;
- analysis of the long-term horizon (strategic level of control) $k = 19$, $k = 30$.

In Fig. 5 to analyze the autocorrelation functions , we shows $W_K^-$ separately

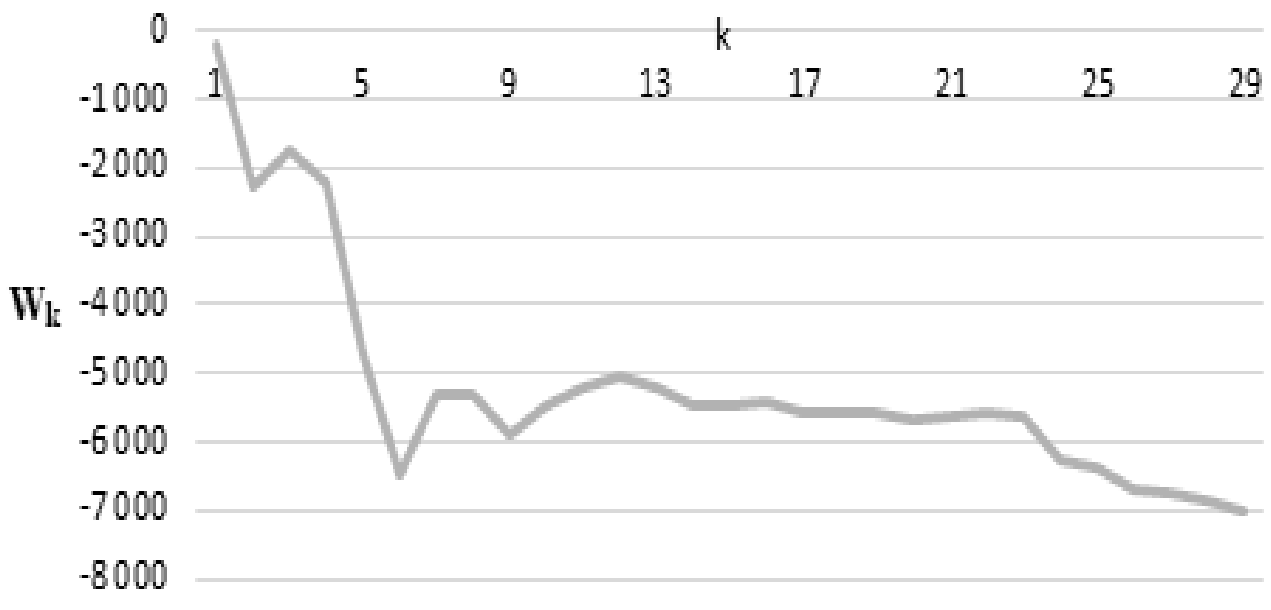

Fig. 5. Three-dimensional modeling of the integral indicator $W_k^-(t)$ of the enterprise with sanctions $k$

In Fig 5The negative curve $W_k^-(t)$ does not repeat the curve $W_k$ in figure 3. It is obvious that the some of the factors giving negative values of the autocorrelation function is the seasonality of the performance of some functions $x^i$, therefore the value $k$ is revised for:

- analysis of the short-term horizon (operational level of control) $k = 2$;
- analysis of the medium-term horizon (tactical level of control) $k = 6$;
- analysis of the long-term horizon (strategic level of control) $k = 24$, $k = 30$.

## VI. The Discussion of the Results

If control $U$ is set for all parameters of the production system, i.e. the dimension of the control vector coincides with the dimension of the dynamic system $X$, then we get the resource distribution system

$$\begin{cases} X^1 = x_1^1 + x_2^1 + ... + x_j^1 + ... x_n^1 + Y^1, \\ X^2 = x_1^2 + x_2^2 + ... + x_j^2 + ... x_n^2 + Y^2, \\ ............................................... \\ X^i = x_1^i + x_2^i + ... + x_j^i + ... x_n^i + Y^i \\ ............................................... \\ X^n = x_1^n + x_2^n + ... + x_j^n + ... x_n^n + Y^n, \end{cases}$$

where

$X^i = [i = 1, 2, ..., n]$ - resource turnover rate;

$Y^i = [i = 1, 2, ..., n]$ - resource turnover at the end node;

$x_j^i = [i = 1, 2, ..., n, j = 1, 2, ..., n]$ - the intensity of the turnover of resources from $i$ the node to create the resource $j$ node.

The controller $a_n^n$ indicates the additional resource consumption for the interaction of one node with another, which forms individual signs of the system from systems of this type. Then we have

$$X^i = \sum_{j=1}^{n} a_j^i X^j + Y^i, i = 1, 2, ..., n.$$

In matrix form, the system has the form

$$X = AX + Y,$$

where

$X = (X^1, X^2, ..., X^n)^T$ – resource intensity vector;

$Y = \left(Y^1, Y^2, ..., Y^n\right)^T$ – vector of intensity of the total resource attracted by;

$$A = \begin{pmatrix} a_1^1 a_2^1 ... a_j^1 ... a_n^1 \\ a_1^2 a_2^2 ... a_j^2 ... a_n^2 \\ a_1^i a_2^i ... a_j^i ... a_n^i \\ a_1^n a_2^n ... a_j^n ... a_n^n \end{pmatrix} = (a_j^i), i, j = \overline{1, n}, \text{ – regulators.}$$

Now it is possible to solve the observation problem from the matrix representation (*E*-*A*)*X*=*Y* (*E* is the identity matrix) and the planning (control) problem $(E\text{-}A)^{-1}Y=X$, which can be solved by optimization methods. Using the proposed (1), we can proceed to the classical form of the control vector of the production system and other tasks: a) set the optimal control solution through the Kalman filter [19], b) go to the classical problem of forecasting the state of a multidimensional dynamic system with an analysis of local functions at each step *t*, c) go to the formal solution of the Cobb-Douglas function for predicting the multidimensional dynamic system by the recurrent equation of state [20], d) go to the solution of Leontief balance models [21], e) solve classical mathematical models of optimal controlled processes and their sufficient optimality conditions [1] and e) apply other optimization methods, discrete-time approaches.

## VII. Conclusion

The method of integral indicators is expanded with the values of the autocorrelation function (6, 7). In Fig. 1 the values of the autocorrelation function for controlling the production system are determined and 2). As a characteristic of developing processes, positive values of the autocorrelation function are grouped in Fig. 4. In Fig. 5 negative values of the autocorrelation function were selected as processes that tend to stop running (the effect of seasonality). According to the developed characteristics, with the help of the regulator, a short-term, medium-term, long-term planning horizon is determined for the operational, tactical, strategic levels of the production system control (in Fig. 3, 4 and 5). With a decrease in value, the amount of information for processing in making managerial decisions decreases, which allows you to timely respond to changes in the external environment during the operational control of the production system.

The tasks set at the beginning of the work are completed:

- formalize the production system as a dynamic system with discrete time;
- develop a digital copy of a dynamic system (production system);
- calculate the autocorrelation function of a space that displays the state of a digital dynamic system (production system).

The purpose of this work, the addition of the integral indicators method to the values of the autocorrelation function to determine the planning horizon of the state of the production system under the influence of the external environment, has been achieved.